\documentclass[12pt,paper,draft]{amsart}

\usepackage{amssymb,latexsym}
\usepackage[matrix,arrow]{xy}
\usepackage[usenames]{color}
\usepackage{slashed}
\usepackage{parskip}
 \numberwithin{equation}{section}
\newtheorem{lemma}{Lemma}[section]

\def\proof{{\indent\it Proof.}}
\def\endproof{\hfill\hbox{$\sqcup$}\llap{\hbox{$\sqcap$}}\medskip}

\begin{document}
\title[ Lower bound estimates for eigenvalues]
{ Lower bound estimates for eigenvalues \\ of  the Laplacian*}
\author [Q. -M. Cheng and X. Qi]{ Qing-Ming Cheng  and Xuerong Qi}
\address{Qing-Ming Cheng\\  Department of Applied Mathematics,
Faculty of Sciences, Fukuoka University, Fukuoka 814-0180, Japan,
cheng@fukuoka-u.ac.jp}
\address{Xuerong Qi\\  Department of Mathematics,
Zhengzhou University, Zhengzhou 450052, P.R. China,
qixuerong609@gmail.com} \subjclass{}
\renewcommand{\thefootnote}{\fnsymbol{footnote}}



\footnotetext{* Research partially supported by a Grant-in-Aid for
Scientific Research from JSPS.}

\begin{abstract}

For an $n$-dimensional polytope $\Omega$  in $\mathbb{R}^{n}$, we
study lower bounds for eigenvalues of the Dirichlet eigenvalue
problem of the Laplacian. In the asymptotic formula on the average
of the first $k$ eigenvalues, Li and Yau \cite{[LY]} obtained the
first term with the order $k^{\frac2n}$, which is optimal. The next
landmark goal is to give the second term with the order
$k^{\frac1n}$ in the asymptotic formula. For this purpose,
Kova\v{r}\'{\i}k, Vugalter and Weidl \cite{[KVW]} have made an
important breakthrough in the case of dimension $2$. It is our
purpose to study the $n$-dimensional case for arbitrary dimension
$n$. We obtain  the second term in the asymptotic sense.
\end{abstract}
\maketitle
\renewcommand{\sectionmark}[1]{}

\section{introduction}

 Let $\Omega \subset \mathbb{R}^{n}$ be a bounded domain with a piecewise smooth boundary
 $\partial \Omega$    in an
 $n$-dimensional Euclidean space $\mathbb{R}^{n}$, $ n\geq 2$. We consider the
 following Dirichlet eigenvalue problem of the Laplacian:
\begin{equation}
\left \{ \aligned \Delta u&=-\lambda u \quad \text{in \ \ $\Omega$}, \\
u&=0 \quad \ \ \ \ \text{on $\partial \Omega$}. \endaligned \right.
\end{equation}
 It is well known that the spectrum of this problem is   real
 and  discrete:
 $$
 0<\lambda_1< \lambda_2\leq \lambda_3\leq \cdots\longrightarrow +\infty,
$$
where each $\lambda_i$ has finite multiplicity which is repeated
according to its multiplicity.

Let $V(\Omega)$ denote the volume of $\Omega$ and let $B_n$ denote
the volume of the unit ball in $\mathbb{R}^n$. One has
 the following Weyl's asymptotic formula
\begin{equation}\lambda_k\sim\frac{4\pi^2}{(B_nV(\Omega))^{\frac{2}{n}}}k^{\frac{2}{n}},
\ \ \ k\rightarrow+\infty.\end{equation} From the above asymptotic
formula, one can obtain
\begin{equation}
\frac{1}{k}\sum_{j=1}^k\lambda_j\sim\frac{n}{n+2}\frac{4\pi^2}{(B_nV(\Omega))^{\frac{2}{n}}}k^{\frac{2}{n}},
\ \ \ k\rightarrow +\infty.
\end{equation}
Furthermore,  P\'olya \cite{[P]} proved that
\begin{equation}\lambda_k\geq\frac{4\pi^2}{(B_nV(\Omega))^{\frac{2}{n}}}k^{\frac{2}{n}},\
\ \ {\rm for} \  k=1,2,\cdots,\end{equation} if $\Omega$ is a tiling
domain in $\mathbb{R}^n$. Moreover,  he proposed  the following

\par \noindent{\bf Conjecture of P\'olya}. {\it
If $\Omega$ is a bounded domain in $\mathbb{R}^n$, then the $k$-th
eigenvalue $\lambda_k$ of the eigenvalue problem {\rm (1.1)}
satisfies
\begin{equation}\lambda_k\geq\frac{4\pi^2}{(B_nV(\Omega))^{\frac{2}{n}}}k^{\frac{2}{n}},\
\ \ {\rm for}\ k=1,2,\cdots.
\end{equation}}
On the conjecture of P\'olya, much work has been done
(\cite{[B]},\cite{[LY]}, \cite{[L]}). In particular, Li and Yau
\cite{[LY]} proved the following
\begin{equation}
\frac{1}{k}\sum_{j=1}^k\lambda_j\geq\frac{n}{n+2}\frac{4\pi^2}{(B_nV(\Omega))^{\frac{2}{n}}}k^{\frac{2}{n}},\
\ \ {\rm for} \ k=1,2,\cdots.\end{equation}
 The formula (1.3) shows
that the constant in the result (1.6) of Li and Yau can not be
improved.
 From this formula (1.6), one can derive
\begin{equation}\label{LY2}
\lambda_k\geq\frac{n}{n+2}\frac{4\pi^2}{(B_nV(\Omega))^{\frac{2}{n}}}k^{\frac{2}{n}},\
\ \ {\rm for} \ k=1,2,\cdots,\end{equation} which gives a partial
solution for the conjecture of P\'olya with a factor
$\frac{n}{n+2}$. Recently, Melas \cite{[M]} improved the estimate
(1.6) to the following:
\begin{equation}
\frac{1}{k}\sum_{j=1}^k\lambda_j\geq\frac{n}{n+2}\frac{4\pi^2}{(B_nV(\Omega))^{\frac{2}{n}}}k^{\frac{2}{n}}
+M_n\frac{V(\Omega)}{I(\Omega)},\ \ \ {\rm for} \
k=1,2,\cdots,\end{equation} where $M_n$ is a positive constant
depending only on the dimension $n$ and
$$
I(\Omega)=\min\limits_{a\in \mathbb{R}^n}\int_\Omega|x-a|^2dx
$$
 is called {\it  the moment of inertia} of
$\Omega$.

 For the average of the first $k$ eigenvalues, it is important
 to compare its lower bound with the following asymptotical behavior:
\begin{equation}
\frac{1}{k}\sum_{j=1}^k\lambda_j=\frac{n}{n+2}\frac{4\pi^2}{(B_nV(\Omega))^{\frac{2}{n}}}k^{\frac{2}{n}}
+C_n\frac{A(\partial\Omega)  \ \ }{ \
V(\Omega)^{1+\frac{1}{n}}}k^{\frac1n}+o(k^{\frac1n}), \ \  \
k\rightarrow +\infty,
\end{equation}
where $A(\partial\Omega)$ denotes the $(n-1)$-dimensional volume of
$\partial\Omega$ and $C_n$ is a positive constant depending only on
the dimension $n$. The first term in (1.9) is due to Weyl
\cite{[W]}. In \cite{[Sa]}, the second term in (1.9) was established
under suitable conditions on $\Omega$. Since the first asymptotical
term is optimal,  the next landmark goal on its lower bound estimate
is to obtain the second asymptotical term with the order of
$k^{\frac1n}$. For this purpose, Kova\v{r}\'{i}k, Vugalter and
 Weidl \cite{[KVW]} have made an important  breakthrough for this landmark goal in the
 case of dimension 2.  They have added  a positive term  in the  right hand side
 of (1.6), which is  similar to  the second term of  (1.9) in the asymptotic sense.
The purpose of this paper is to study  the $n$-dimensional case  for
arbitrary dimension $n$.  We also obtain the second term of (1.9) in
the asymptotic sense. For estimates on upper bounds of eigenvalues,
one can see Cheng and Yang \cite{[cy]}.

For an $n$-dimensional polytope $\Omega$ in $\mathbb{R}^n$, we
denote by $p_i, \ i=1,\cdots,m,$ the $i$-th face of $\Omega$. Assume
that $A_i$ is the area of the $i$-th face $p_i$ of $\Omega$. For
each $i=1,\cdots,m$, we choose several non-overlapping
$(n-1)$-dimensional convex subdomains $s_{r_i}$ in the interior of
$p_i$ such that the area of $\bigcup s_{r_i}$ is greater than or
equal to one third of $A_i$ and the distance $d_i$ between $\bigcup
s_{r_i}$ and $\partial\Omega\setminus p_i$ is greater than 0. Define
the function $\Theta:\mathbb{R}\rightarrow\mathbb{R}$ by
$\Theta(t)=0$ if $t\leq0$ and $\Theta(t)=1$ if $t>0$. Then we prove
the following

\par\noindent{\bf Theorem 1}. {\it Let $\Omega$ be an
$n$-dimensional polytope in $\mathbb{R}^n$. Then, for any positive
integer $k$, we have
$$\frac{1}{k}\sum_{j=1}^{k} \lambda_j
\geq\frac{n}{n+2}\frac{4\pi^2}{(B_nV(\Omega))^{\frac2n}}k^{\frac{2}{n}}
+\frac{3^{-4}2^{3-n}\pi^2}{(n+2)B_n^{\frac2n}}
\frac{A(\partial\Omega) \ \ }{ \ V(\Omega)^{1+\frac
2n}}\biggl(\frac{V(\Omega)\lambda_k}{\alpha_1}\biggl)^{-n\varepsilon(k)}k^{\frac{2}{n}}\lambda_k^{-\frac{1}{2}}\Theta
(\lambda_k-\lambda_0), $$ where
$$\aligned\varepsilon(k)&=\left[\sqrt{\frac{{\rm
log}_{2}\big((V(\Omega)/\alpha_1)^{n-1} \lambda_k^{\frac
n2}\big)}{n+12}}~\right]^{-1} ,\quad
  \alpha_1=\sqrt{\frac{3}{B_n}\left(\frac{4n\pi^2}{n+2}\right)^\frac
n2},\\
\lambda_0&={\rm max}\left\{\frac{4n}{\underset{i}{\rm
min}\{d_i^2\}}, \
\left(\frac{\alpha_1}{V(\Omega)}\right)^{\frac{2}{n}}, \
2^{\frac{2(n+12)}{n}}\left(\frac{\alpha_1}{V(\Omega)}\right)^{\frac{2(n-1)}{n}},
\ \left(\frac{12}{\underset{i}{\rm min}
\{A_i\}}\right)^{\frac{2}{n-1}}\right\}.\endaligned$$}

\par\noindent{\bf Remark 1}. Notice that
$\varepsilon(k)\rightarrow0$ and
$\lambda_k\sim\frac{4\pi^2}{(B_nV(\Omega))^{\frac{2}{n}}}k^{\frac{2}{n}}$
as $k\rightarrow+\infty$. It shows that the second term on the right
hand side of the inequality in Theorem 1 is very similar to the
second term in the asymptotic (1.9) when $k$ is large enough.
Combining this with the inequality (1.7), we immediately obtain the
following

\par\noindent{\bf Corollary 1}. {\it Let $\Omega$ be an
$n$-dimensional polytope in $\mathbb{R}^n$. Then, there exists a
positive integer $N$, such that, for all $k\geq N$,
$$\frac{1}{k}\sum_{j=1}^{k} \lambda_j
\geq\frac{n}{n+2}\frac{4\pi^2}{(B_nV(\Omega))^{\frac2n}}k^{\frac{2}{n}}
+\frac{\pi }{3^42^{n-1}(n+2)B_n^{\frac1n}} \frac{A(\partial\Omega) \
\ }{ \ V(\Omega)^{1+\frac 1n}}k^{\frac{1}{n}-2\varepsilon(k)},
$$
where
$$\varepsilon(k)=\left[\sqrt{\frac{1}{n+12}{\rm
log}_{2}\left(\bigg(\frac{V(\Omega)}{\alpha_1}\bigg)^{n-1}
\bigg(\frac{4n\pi^2}{n+2}\bigg)^{\frac
n2}\frac{k}{B_nV(\Omega)}\right)}~\right]^{-1},\quad
  \alpha_1=\sqrt{\frac{3}{B_n}\left(\frac{4n\pi^2}{n+2}\right)^\frac
n2}.$$}


\section{Proof of Main theorem}

For an $n$-dimensional polytope $\Omega$ in $\mathbb{R}^n$, let
$u_j$ be a normalized eigenfunction corresponding to the $j$-th
eigenvalue $\lambda_j$, i.e. $u_j$ satisfies
\begin{equation}
\left \{ \aligned &\Delta u_j=-\lambda_j u_j \qquad \text{in \ \ $\Omega$}, \\
& ~ ~ \ ~ \  u_j=0 \qquad\quad \ \ \ \text{on $\partial \Omega$},\\
&\int_\Omega u_ju_k=\delta_{jk}, \quad  \ \forall  \ j, k.
\endaligned \right.
\end{equation}
 Then $\{u_j\}_{j=1}^{+\infty}$ forms an orthonormal basis of
$L^2(\Omega)$. We consider the function $\varphi_j$ given by
\begin{equation*}
\varphi_{j}(x)= \left \{ \aligned u_j(x) \ \ \ \  , \ \ \
\ & x\in \Omega,\\
 0 \ \ \ \ \ \  , \ \ \ \ & x\in\mathbb{R}^n\setminus\Omega.
\endaligned \right.
\end{equation*}
Denote by $\widehat{\varphi}_j$ the Fourier transform of
$\varphi_j$. For any $\xi\in\mathbb{R}^n$, we have
$$\widehat{\varphi}_j(\xi)=(2\pi)^{-\frac n2}\int_{\mathbb{R}^n}\varphi_j(x)e^{\sqrt{-1}<\xi, x>}dx
=(2\pi)^{-\frac n2}\int_{\Omega}u_j(x)e^{\sqrt{-1}<\xi, x>}dx.$$
Take $\lambda$ large enough such that $\lambda$ satisfies the
following two conditions:\\
(i) for each $i=1,\cdots,m$, there exist some non-overlapping
$(n-1)$-dimensional cubes $t_{l_i}$ with the side
$\frac1{\sqrt{\lambda}}$ on $\bigcup s_{r_i}$, whose
total area is greater than or equal to $\frac16 A_i$;\\
(ii) $\frac1{\sqrt{\lambda}}\leq\frac{1}{2\sqrt{n}}\underset{i}{\rm
min}\{d_i\}\leq\frac{d_i}{2\sqrt{n}}.$ In this case, we make sure
that the $n$-dimensional rectangles
$T_{l_i}=[0,\frac{1}{2\sqrt{\lambda}}]\times t_{l_i}$ lie inside
$\Omega$ and they do not overlap each other.

We define a function $F_\lambda$ by
$$F_\lambda(\xi)=\sum_{\lambda_j\leq\lambda}\bigl|\widehat{\varphi}_j(\xi)\bigl|^2.$$
By Parseval's identity, we have
\begin{equation}
\int_{\mathbb{R}^n}F_\lambda(\xi)d\xi=\sum_{\lambda_j\leq\lambda}\int_{\mathbb{R}^n}
|\widehat{\varphi}_j(\xi)|^2d\xi=\sum_{\lambda_j\leq\lambda}\int_{\mathbb{R}^n}
\varphi_j^2(x)dx=\sum_{\lambda_j\leq\lambda}\int_{\Omega}
u_j^2(x)dx=N(\lambda),
\end{equation}
where $N(\lambda)$ is the number of eigenvalues
$\lambda_j\leq\lambda$.  Furthermore, we deduce from integration by
parts and Parseval's identity that
\begin{equation}
\int_{\mathbb{R}^n}|\xi|^{2}F_\lambda(\xi)d\xi=\sum_{\lambda_j\leq\lambda}\int_{\mathbb{R}^n}|\xi|^{2}|\widehat{\varphi}_j(\xi)|^2d\xi
=-\sum_{\lambda_j\leq\lambda}\int_\Omega u_j\Delta
u_jdx=\sum_{\lambda_j\leq\lambda}\lambda_j.
\end{equation}

For each fixed $\xi\in\mathbb{R}^n$, since $e^{\sqrt{-1}<\xi,x>}$
belongs to $L^2(\Omega)$, it follows that
$$e^{\sqrt{-1}<\xi,x>}=\sum_{j=1}^{\infty}c_j(\xi)u_j, \ \ {\rm where} \ \ c_j(\xi)=\int_\Omega u_j(x)e^{\sqrt{-1}<\xi,x>}dx.$$
Let
$$u(\xi,x)=\sum_{\lambda_j\leq\lambda}c_j(\xi)u_j(x).$$
Then we have
\begin{equation}\biggl\|u-e^{\sqrt{-1}<\xi,x>}\biggl\|^2_{L^2(\Omega)}=
V-(2\pi)^nF_\lambda(\xi). \end{equation}

To prove Theorem 1, we will need the following lemma.

\begin{lemma}{\rm(\cite{[LY]})}
If $F$ is a real-valued function defined on $\mathbb{R}^n$ with
$0\leq F\leq M_1$, and
$$\int_{\mathbb{R}^n}|\xi|^2F(\xi)d\xi\leq M_2,$$
then
$$\int_{\mathbb{R}^n}F(\xi)d\xi\leq\bigg(\frac{n+2}{n}\bigg)^{\frac{n}{n+2}}M_2^{\frac{n}{n+2}}
\big(M_1B_n\big)^{\frac{2}{n+2}}.$$
\end{lemma}

Next we need to establish the estimate for $F_\lambda(\xi)$.  For
each $l_i$, we choose a local
 coordinate system $(x_1,\cdots,x_n)$ such that $t_{l_i}=\bigl[-\frac{1}{2\sqrt{\lambda}},\frac{1}{2\sqrt{\lambda}}\bigl]^{n-1}$
 and $\frac{\partial}{\partial x_1}$ is the inward unit normal vector field on
 $ t_{l_i}$.  To
derive the upper bound of $F_\lambda(\xi)$, we prepare the following
lemmas:
\begin{lemma}
 For any positive integer $p$, we have
$$\biggl\|\frac{\partial^{p}u}{\partial x_1^{p}}\biggl\|^{2}_{L^{2}(T_{l_i})}
\leq \biggl(\frac{n+2}{4n\pi^2}\biggl)^{\frac{n}{2}}B_nV^{2}
D_{p-1}\lambda^{p+\frac{n}{2}},$$ where the sequence $D_{q}$ is
defined by
$$\aligned D_0&=1, \qquad D_1=3\bigl(1+44^2 n^{2}
p^{4}+4\cdot 5^2n p^{2}\bigl),\\
 D_q&=3\bigl(1+44^2 n^{2}
p^{4}\bigl)D_{q-2}+\bigl(12\cdot 5^2n p^{2}\bigl)D_{q-1}, \ \
q=2,3,\cdots.\endaligned$$
\end{lemma}

\proof \ \   For $p\geq1$ and $0\leq q\leq p-1$, we define functions
$g$ and $v_{q,p}$ by the following
$$
g(x)=1-6x^{4}+8x^6-3x^8, \ \ \  0\leq x\leq 1, \ \ \ \ \ \ \ \ \ \ \ \ \ \ \ \ \ \
$$
\begin{equation*}
v_{q,p}(t)=
\left \{ \aligned 1 \ \ \ \ \ \ \ \ \ \ , \ \ \ \ & 0 \leq t \leq\frac{2p-q}{2p}, \\
g(2pt-2p+q) \ \ , \ \ \ \ & \frac{2p-q}{2p}\leq t\leq \frac{2p-q+1}{2p},\\
0 \ \ \ \ \ \ \ \ , \ \ \ \ & \frac{2p-q+1}{2p}< t,
\endaligned \right.
\end{equation*}
with $v_{q,p}(-t)=v_{q,p}(t)$ for $t<0$. From the definition of $g$,
it follows that
\begin{equation}\label{g}
|g(x)|\leq 1, \qquad |g'(x)|<\frac{5}{2}, \qquad |g''(x)|< 11.
\end{equation}
By the definition of $v_{q,p}$ and \eqref{g}, we get
\begin{equation}\label{v}
|v_{q,p}(t)|\leq 1, \ \ \ \ |v_{q,p}'(t)|< 5p, \ \ \ \
|v_{q,p}''(t)|< 44 p^{2}.
\end{equation}
 Next we define
$$
W_{q,p,\lambda}(x_1,\cdots,x_n)= v_{q,p}(\sqrt{\lambda}~
x_1)v_{q,p}(\sqrt{\lambda}~ x_2)\cdots v_{q,p}(\sqrt{\lambda}~ x_n),
\ \ \ \ (x_1,\cdots, x_n)\in \mathbb{R}^n, $$ and set
$$\omega_q={\rm supp} W_{q,p,\lambda}.  $$ Then we have
\begin{equation}\label{W}
\big|W_{q,p,\lambda}\big|\leq1, \ \ \ \ \big|\nabla
W_{q,p,\lambda}\big|<5n^{\frac12}\lambda^{\frac12}p, \ \ \ \
\big|\Delta W_{q,p,\lambda}\big|<44n \lambda p^2,
\end{equation}
and \begin{equation}\label{om}
 T_{l_i}\subset\omega_q, \quad
\omega_{q+1}\subset \omega_{q}, \quad  W_{q,p,\lambda}\equiv 1 \
{\rm on} \ \omega_{q+1}, \quad \nabla W_{q,p,\lambda}=0 \ {\rm on} \
\partial\omega_q.
\end{equation}

We will prove
\begin{equation}\label{in}
 \biggl\|\nabla\bigg(\frac{\partial^{q}\varphi_j}{\partial x_1^{q}}\bigg)\biggl\|^{2}_{L^{2}(\omega_q)}\leq D_{q}\lambda^{q+1},
\end{equation} by induction on $q$ for $q=0,1,\cdots,p-1$.
For $q=0$,  since $D_0=1$, it follows that
\begin{equation}\label{q=0}
\big\|\nabla \varphi_j\big\|^2_{L^2(\omega_0)}\leq
\int_\Omega|\nabla u_j|^2=-\int_\Omega u_j\Delta
u_j=\lambda_j\leq\lambda= D_{0}\lambda.
\end{equation}
Using the property \eqref{om} of $W_{q,p,\lambda}$ and integration
by parts, we get
\begin{equation}\label{q=1}\aligned
&\biggl\|\Delta\big(\varphi_jW_{0,p,\lambda}\big)\biggl\|^2_{L^2(\omega_{0})}\\
=\sum_{k=1}^{n}&\biggl\|\frac{\partial^2}{\partial
x_k^2}\big(\varphi_jW_{0,p,\lambda}\big)\biggl\|^2_{L^2(\omega_{0})}
 +2\sum_{k<l}\int_{\omega_{0}}\frac{\partial^2}{\partial
x_k^2}\big(\varphi_jW_{0,p,\lambda}\big)\frac{\partial^2}{\partial
x_l^2}\big(\varphi_jW_{0,p,\lambda}\big)\\
=\sum_{k=1}^{n}&\biggl\|\frac{\partial^2}{\partial
x_k^2}\big(\varphi_jW_{0,p,\lambda}\big)\biggl\|^2_{L^2(\omega_{0})}+2\sum_{k<l}\biggl\|\frac{\partial^2}{\partial
x_k
\partial x_l}\big(\varphi_jW_{0,p,\lambda}\big)\biggl\|_{L^2(\omega_{0})}^2.
\endaligned\end{equation}

From \eqref{W}, \eqref{om}, \eqref{q=0} and \eqref{q=1}, it follows
that
$$
\aligned \biggl\|\nabla\bigg(\frac{\partial\varphi_j}{\partial
x_1}\bigg)\biggl\|^{2}_{L^{2}(\omega_{1})} &=\sum_{k=1}^{n}
\biggl\|\frac{\partial^{2}}{\partial x_1\partial
x_{k}}\bigl(\varphi_jW_{0,p,\lambda}\bigl)\biggl\|^{2}_{L^{2}(\omega_{1})}\\
 &\leq  \sum_{k=1}^{n} \biggl\|\frac{\partial^{2}}{\partial
x_1\partial x_{k}}\bigl(\varphi_jW_{0,p,\lambda}\bigl)\biggl\|^{2}_{L^{2}(\omega_{0})}\\
 &\leq
\biggl\|\Delta\big(\varphi_jW_{0,p,\lambda}\big)\biggl\|^2_{L^2(\omega_{0})}\\
 &=\biggl\|-\lambda_j\varphi_jW_{0,p,\lambda}+\varphi_j\Delta
W_{0,p,\lambda}+2\nabla\varphi_j\nabla
W_{0,p,\lambda}\biggl\|^2_{L^2(\omega_{0})}\\
 &\leq
 3\biggl\{\lambda^2\big(1+44^2 n^2 p^4
\big)\big\|\varphi_j\big\|^2_{L^2(\omega_{0})}+4\cdot
5^2p^2n\lambda\big\|\nabla\varphi_j\big\|^2_{L^2(\omega_{0})}\biggl\}\\
&\leq 3\biggl\{\lambda^2\bigl(1+44^2 n^2 p^4 \bigl)+4\cdot
5^2p^2n\lambda^2\biggl\} =D_{1}\lambda^{2}.
\endaligned
$$
Hence \eqref{in} holds for $q=0$ and $q=1$.  Now assume that
\eqref{in} holds some $q-1$ and $q$. We will show that it holds for
$q+1$ as well. Notice that
$$W_{q,p,\lambda}=\nabla W_{q,p,\lambda}=0 \ \ \ \  {\rm on} \ \ \partial\omega_{q},$$ then
\begin{equation}\label{q=q}\aligned
&\biggl\|\Delta\biggl(\frac{\partial^{q}\varphi_j}{\partial
x_1^{q}}W_{q,p,\lambda}\biggl)\biggl\|^2_{L^2(\omega_{q})}\\
=\sum_{k=1}^{n}&\biggl\|\frac{\partial^2}{\partial
x_k^2}\biggl(\frac{\partial^{q}\varphi_j}{\partial
x_1^{q}}W_{q,p,\lambda}\biggl)\biggl\|^2_{L^2(\omega_{q})}
 +2\sum_{k<l}\biggl\|\frac{\partial^2}{\partial x_k
\partial x_l}\biggl(\frac{\partial^{q}\varphi_j}{\partial
x_1^{q}}W_{q,p,\lambda}\biggl)\biggl\|_{L^2(\omega_{q})}^2.
\endaligned\end{equation}
From \eqref{W}, \eqref{om} and \eqref{q=q}, it follows that
$$\aligned & \ \ \ \
\biggl\|\nabla\bigg(\frac{\partial^{q+1}\varphi_j}{\partial
x_1^{q+1}}\bigg)\biggl\|^{2}_{L^{2}(\omega_{q+1})}\\
&=\sum_{k=1}^{n}\biggl\|\frac{\partial^2}{\partial x_1\partial
x_k}\biggl(\frac{\partial^{q}\varphi_j}{\partial
x_1^{q}}W_{q,p,\lambda}\biggl)\biggl\|^2_{L^2(\omega_{q+1})}\\
&\leq \sum_{k=1}^{n}\biggl\|\frac{\partial^2}{\partial x_1\partial
x_k}\biggl(\frac{\partial^{q}\varphi_j}{\partial
x_1^{q}}W_{q,p,\lambda}\biggl)\biggl\|^2_{L^2(\omega_{q})}\\
& \leq \biggl\|\Delta\biggl(\frac{\partial^{q}\varphi_j}{\partial
x_1^{q}}W_{q,p,\lambda}\biggl)\biggl\|^2_{L^2(\omega_{q})}\\
 &=\biggl\|-\lambda_j\biggl(\frac{\partial^{q}\varphi_j}{\partial
x_1^{q}}\biggl)W_{q,p,\lambda}+\biggl(\frac{\partial^{q}\varphi_j}{\partial
x_1^{q}}\biggl)\Delta
W_{q,p,\lambda}+2\nabla\biggl(\frac{\partial^{q}\varphi_j}{\partial
x_1^{q}}\biggl)\nabla
W_{q,p,\lambda}\biggl\|^2_{L^2(\omega_{q})}\\
 &\leq 3\biggl\{\lambda^2\bigl(1+44^2 n^2 p^4
\bigl)\biggl\|\frac{\partial^{q}\varphi_j}{\partial
x_1^{q}}\biggl\|^2_{L^2(\omega_{q})}+4\cdot
5^2p^2n\lambda\biggl\|\nabla\frac{\partial^{q}\varphi_j}{\partial
x_1^{q}}\biggl\|^2_{L^2(\omega_{q})}\biggl\}\\
&\leq 3\biggl\{\lambda^2\bigl(1+44^2 n^2 p^4
\bigl)\biggl\|\nabla\bigg(\frac{\partial^{q-1}\varphi_j}{\partial
x_1^{q-1}}\bigg)\biggl\|^2_{L^2(\omega_{q-1})}+4\cdot
5^2p^2n\lambda\biggl\|\nabla\frac{\partial^{q}\varphi_j}{\partial
x_1^{q}}\biggl\|^2_{L^2(\omega_{q})}\biggl\}\\
&\leq \lambda^{q+2}\biggl\{3\bigl(1+44^2 n^2
p^4\bigl)D_{q-1}+12\cdot 5^2np^2D_{q}\biggl\}=D_{q+1}\lambda^{q+2}.
\endaligned$$
This implies that \eqref{in} holds for $q+1$. Therefore, \eqref{in}
holds for any integer $0\leq q\leq p-1$. Taking $q=p-1$ in
\eqref{in}, we can obtain
\begin{equation}
\biggl\|\frac{\partial^{p}u_j}{\partial
x_1^{p}}\biggl\|^{2}_{L^{2}(T_{l_i})}\leq
D_{p-1}\lambda^{p}.\end{equation} Using the Cauchy-Schwarz
inequality, we have
$$\aligned
\biggl\|\frac{\partial^{p}u}{\partial
x_1^{p}}\biggl\|^{2}_{L^{2}(T_{l_i})}&=\biggl\|\sum_{\lambda_j\leq
\lambda}c_j\frac{\partial^{p}u_j}{\partial
x_1^{p}}\biggl\|^{2}_{L^{2}(T_{l_i})}=\int_{T_{l_i}}\biggl
|\sum_{\lambda_j\leq\lambda}c_j\frac{\partial^{p}u_j}{\partial
x_1^{p}}\biggl |^2\\
&\leq
\int_{T_{l_i}}\sum_{\lambda_j\leq\lambda}|c_j|^2\sum_{\lambda_j\leq\lambda}\biggl
|\frac{\partial^{p}u_j}{\partial x_1^{p}}\biggl |^2\leq V
\sum_{\lambda_j\leq\lambda}\biggl \|\frac{\partial^{p}u_j}{\partial
x_1^{p}}\biggl \|^2_{L^2(T_{l_i})}\leq V\cdot
N(\lambda)D_{p-1}\lambda^{p}.\endaligned$$
 Using the lower bound of
$\lambda_j$ given in \eqref{LY2}, we find out that
\begin{equation*}
N(\lambda)\leq
B_nV\biggl(\frac{n+2}{4n\pi^2}\biggl)^{\frac{n}{2}}\lambda^{\frac{n}{2}}.\end{equation*}
 Then
$$ \ \ \ \ \ \ \ \ \ \ \ \ \ \ \ \ \ \ \ \ \ \ \ \ \ \ \ \ \ \ \ \ \
\ \ \ \ \ \ \biggl\|\frac{\partial^{p}u}{\partial
x_1^{p}}\biggl\|^{2}_{L^{2}(T_{l_i})}\leq
\biggl(\frac{n+2}{4n\pi^2}\biggl)^{\frac{n}{2}}B_nV^{2}
D_{p-1}\lambda^{p+\frac{n}{2}}. \ \ \ \ \ \ \ \ \ \ \ \ \ \ \ \ \ \
\ \ \ \ \ \ \ \ \ \
 \endproof$$

\begin{lemma}
Let $p$ be a positive integer and let $f\in
C^{p}\bigl[0,\frac{1}{2\sqrt{\lambda}}\bigl]$ be a real-valued
function. If $f^{(p)}$ is not identically zero, then one of the
following inequalities holds true:
$$\aligned
{\rm max}|f'| &\leq 2^{p-1}\big({\rm
max}|f^{(p)}|\big)^{\frac{1}{p}}\big({\rm
max}|f|\big)^{1-\frac{1}{p}},\\
{\rm max}\big|f'\big| &< 4^{p+1}\lambda^{\frac{1}{2}}{\rm
max}|f|.\endaligned$$
\end{lemma}

\proof \ \ Let $m_q={\rm max}|f^{(q)}|$, $q\in\{0,1,\cdots,p\}$. If
$p=1$, then the conclusion is obvious. Next we assume that $p\geq
2$. For any fixed $0\leq q\leq p-2$, we let
$m_{q+1}=|f^{(q+1)}(t_0)|$, $t_0 \in
\bigl[0,\frac{1}{2\sqrt{\lambda}}\bigl]$. We discuss separately the
following cases:

Case 1:  $t_0<\frac{1}{4\sqrt{\lambda}}$. In this case, if
$\frac{m_{q+1}}{m_{q+2}}\leq\frac{1}{4\sqrt{\lambda}}$, then it
follows from Taylor's formula that
$$\aligned
m_q\geq \bigg|f^{(q)}\bigg(t_0+\frac{m_{q+1}}{m_{q+2}}\bigg)\bigg|
&=\biggl|f^{(q)}(t_0)+f^{(q+1)}(t_0)\cdot\frac{m_{q+1}}{m_{q+2}}+\frac{1}{2}f^{(q+2)}
(\xi)\biggl(\frac{m_{q+1}}{m_{q+2}}\biggl)^2\biggl|\\
&\geq\bigg|f^{(q+1)}(t_0)\cdot\frac{m_{q+1}}{m_{q+2}}\bigg|-\bigg|f^{(q)}(t_0)\bigg|-\bigg|\frac{1}{2}f^{(q+2)}
(\xi)\biggl(\frac{m_{q+1}}{m_{q+2}}\biggl)^2\biggl|\\
 &\geq
-m_q+\frac{m_{q+1}^2}{m_{q+2}}-\frac{m_{q+2}}{2}\biggl(\frac{m_{q+1}}{m_{q+2}}\biggl)^2\\
&=-m_q+\frac{m_{q+1}^2}{2m_{q+2}},
\endaligned$$
which implies
\begin{equation}\label{c11}
\frac{m_{q+2}}{m_{q+1}}\geq \frac{1}{4}\frac{m_{q+1}}{m_q}.
\end{equation}
If $\frac{m_{q+1}}{m_{q+2}}>\frac{1}{4\sqrt{\lambda}}$,
then we have
$$\aligned
m_q\geq\bigg|
f^{(q)}\biggl(t_0+\frac{1}{4\sqrt{\lambda}}\biggl)\bigg|&=
\biggl|f^{(q)}(t_0)+f^{(q+1)}(t_0)\cdot\frac{1}{4\sqrt{\lambda}}+\frac{1}{2}f^{(q+2)}(\xi)\biggl(\frac{1}{4\sqrt{\lambda}}\biggl)^2\biggl|\\
&\geq
-m_q+\frac{m_{q+1}}{4\sqrt{\lambda}}-\frac{m_{q+2}}{2}\biggl(\frac{1}{4\sqrt{\lambda}}\biggl)^2\\
&>
-m_q+\frac{m_{q+1}}{4\sqrt{\lambda}}-\frac{4\sqrt{\lambda}~m_{q+1}}{2}\biggl(\frac{1}{4\sqrt{\lambda}}\biggl)^2,
\endaligned$$
which implies
\begin{equation}\label{c12}
 \frac{m_{q+1}}{m_q}<
16\sqrt{\lambda}.\end{equation}
Case 2:
$t_0\geq\frac{1}{4\sqrt{\lambda}}.$ In this case, if
$\frac{m_{q+1}}{m_{q+2}}\leq \frac{1}{4\sqrt{\lambda}}$,
 then we have
$$\aligned
m_q\geq \bigg|
f^{(q)}(t_0-\frac{m_{q+1}}{m_{q+2}})\bigg|&=\biggl|f^{(q)}(t_0)
+f^{(q+1)}(t_0)\biggl(-\frac{m_{q+1}}{m_{q+2}}\biggl)+\frac{1}{2}f^{(q+2)}(\xi)\biggl(-\frac{m_{q+1}}{m_{q+2}}\biggl)^2\biggl|\\
&\geq
-m_q+\frac{m_{q+1}^2}{m_{q+2}}-\frac{m_{q+2}}{2}\biggl(\frac{m_{q+1}}{m_{q+2}}\biggl)^2\\
&=-m_q+\frac{m_{q+1}^2}{2m_{q+2}},
\endaligned$$
which implies \eqref{c11}. If $\frac{m_{q+1}}{m_{q+2}}>
\frac{1}{4\sqrt{\lambda}}$,
 then we have
$$\aligned
m_q\geq \bigg|
f^{(q)}(t_0-\frac{1}{4\sqrt{\lambda}})\bigg|&=\biggl|f^{(q)}(t_0)+f^{(q+1)}(t_0)\biggl(-\frac{1}{4\sqrt{\lambda}}\biggl)
+\frac{1}{2}f^{(q+2)}(\xi)\biggl(-\frac{1}{4\sqrt{\lambda}}\biggl)^2\biggl|\\
&\geq
-m_q+\frac{m_{q+1}}{4\sqrt{\lambda}}-\frac{m_{q+2}}{2}\biggl(\frac{1}{4\sqrt{\lambda}}\biggl)^2\\
&>
-m_q+\frac{m_{q+1}}{4\sqrt{\lambda}}-\frac{4\sqrt{\lambda}~m_{q+1}}{2}\biggl(\frac{1}{4\sqrt{\lambda}}\biggl)^2,
\endaligned$$
which implies \eqref{c12}. Therefore, for any fixed $0\leq q\leq
p-2$,  one of the inequalities \eqref{c11} and \eqref{c12} holds
true.

Meanwhile, we note that there are two possibilities. Either for all
$0\leq q\leq p-1$,
\begin{equation}\label{c21}
\frac{m_{q+1}}{m_{q}}\geq 16\sqrt{\lambda},\end{equation} or there
exists $q_0\in \{0,1,\cdots,p-1\},$ such that
\begin{equation}\label{c22}
\forall \ 0\leq q<q_0, \ \ \ \ \frac{m_{q+1}}{m_{q}}\geq
16\sqrt{\lambda}, \ \ \ \ \ \
\frac{m_{q_0+1}}{m_{q_0}}<16\sqrt{\lambda}.\end{equation} If
\eqref{c21} holds, then we apply \eqref{c11} to get
$$\aligned
\frac{m_p}{m_0}&=\frac{m_p}{m_{p-1}}\cdot\frac{m_{p-1}}{m_{p-2}}\cdots\frac{m_2}{m_1}\cdot\frac{m_1}{m_0}\\
&\geq
\frac{1}{4}\frac{m_{p-1}}{m_{p-2}}\cdot\frac{1}{4}\frac{m_{p-2}}{m_{p-3}}\cdots\frac{1}{4}\frac{m_{1}}{m_{0}}\cdot\frac{m_1}{m_{0}}\\
&=\bigg(\frac{1}{4}\bigg)^{p-1}\cdot\frac{m_{p-1}}{m_0}\cdot\frac{m_1}{m_0},
\endaligned$$
which implies
$$\frac{m_p}{m_0}\geq \biggl(\frac{1}{4}\biggl)^{(p-1)+\cdots +1}\biggl(\frac{m_1}{m_0}\biggl)^{p}=4^{-\frac{p(p-1)}{2}}
\biggl(\frac{m_1}{m_0}\biggl)^{p}.$$ Thus,
$$\frac{m_1}{m_0}\leq 2^{p-1}\biggl(\frac{m_p}{m_0}\biggl)^{\frac{1}{p}}.$$
If \eqref{c22} holds, we apply \eqref{c11} to obtain
$$ \frac{m_1}{m_0}\leq 4\frac{m_2}{m_1}\leq \cdots \leq 4^{q_0}\frac{m_{q_0+1}}{m_{q_0}}
<4^{q_0}\cdot 16\sqrt{\lambda}=4^{q_0+2}\lambda^{\frac12}\leq
4^{p+1}\lambda^{\frac12}. $$ This completes the proof.
\endproof

\begin{lemma}
Let $p$ be a positive integer. Then for any $\xi\in\mathbb{R}^n,$
$$\biggl\|u-e^{\sqrt{-1}<\xi, x>}\biggl\|^2_{L^{2}(T_{l_i})} \geq \frac{1}{9\cdot 2^{n-1}}{\rm
min}\biggl \{2^{-2p-5}\lambda^{-\frac n2}, \ \ 2^{-p-2
}6^{\frac{1}{2p}}(\beta_p^2+\beta_{p+1}^2)^{-\frac{1}{2p}}\lambda^{-\frac
n2-\frac{n}{2p}} \biggl\},$$
where
$$\beta_{q}^2=\biggl(\frac{n+2}{4n\pi^2}\biggl)^{\frac{n}{2}}B_nV^{2}
D_{q-1}, \ \ \ \ q\in \{p,p+1\}.$$
\end{lemma}

\proof \ \ By Lemma 2.2, we have
$$\biggl\|\frac{\partial^{q}u}{\partial x_1^{q}}\biggl\|^{2}_{L^{2}(T_{l_i})}
\leq \beta_{q}^2\lambda^{q+\frac{n}{2}}.$$
 This implies that the measure of the set
$$\biggl\{(x_2,\cdots,x_n)\in\biggl[-\frac{1}{2\sqrt{\lambda}},\frac{1}{2\sqrt{\lambda}}\biggl]^{n-1}: \ \int_0^{\frac{1}{2\sqrt{\lambda}}}
\biggl |\frac{\partial^{q}u}{\partial x_1^{q}}\biggl |^2dx_1\leq
2\beta_q^2\lambda^{q+n-\frac{1}{2}}, \ \ q\in\{p, p+1\}\biggl\}$$ is
obviously at least $(\frac{1}{2\sqrt{\lambda}})^{n-1}$. For such
$(x_2,\cdots,x_n)$, we let $\underset{x_1}{\rm max}
\big|\frac{\partial^{p}u}{\partial
x_1^{p}}\big|=\big|\frac{\partial^{p}u}{\partial
x_1^{p}}(x_1^0)\big|$
 with $x_1^0\in [0,\frac{1}{2\sqrt{\lambda}}]$.
 For any $x_1\in[0,\frac{1}{2\sqrt{\lambda}}]$, we have
$$\left|\frac{\partial^{p}u}{\partial x_1^{p}}(x_1^0)\right|^2
=\biggl|\frac{\partial^{p}u}{\partial
x_1^{p}}(x_1)-\int_{x_1^0}^{x_1}\frac{\partial^{p+1}u}{\partial
x_1^{p+1}}(\tau)d\tau\biggl|^2.$$
 Integrating both sides of the equality with respect to $x_1$ and using Jensen's inequality,
we obtain
$$
\aligned \ \ \ \ \ \
&\frac{1}{2\sqrt{\lambda}}\bigg(\underset{{x_1}}{\rm
max}\biggl|\frac{\partial^{p}u}{\partial x_1^{p}}\biggl|\bigg)^2\\
 =&\int_{0}^{\frac{1}{2\sqrt{\lambda}}} \biggl
|\frac{\partial^{p}u}{\partial
x_1^{p}}-\int_{x_1^0}^{x_1}\frac{\partial^{p+1}u}{\partial
x_1^{p+1}}(\tau)d\tau\biggl |^2dx_1\\
=&\int_{0}^{\frac{1}{2\sqrt{\lambda}}} \bigg
|\frac{\partial^{p}u}{\partial
x_1^{p}}\bigg|^2dx_1+\int_0^{\frac{1}{2\sqrt{\lambda}}}\biggl
|\int_{x_1^0}^{x_1}\frac{\partial^{p+1}u}{\partial
x_1^{p+1}}(\tau)d\tau\biggl |^2 dx_1\\
& -2{\rm
Re}\int_0^{\frac{1}{2\sqrt{\lambda}}}\biggl[\frac{\partial^{p}u}{\partial
x_1^{p}} \overline{\int_{x_1^0}^{x_1}\frac{\partial^{p+1}u}{\partial
x_1^{p+1}}(\tau)d\tau}~\biggl ]dx_1\\
\endaligned$$
$$\aligned
\leq &\int_{0}^{\frac{1}{2\sqrt{\lambda}}}
\bigg|\frac{\partial^{p}u}{\partial
x_1^{p}}\bigg|^2dx_1+\int_0^{\frac{1}{2\sqrt{\lambda}}}\biggl
|\int_{x_1^0}^{x_1}\frac{\partial^{p+1}u}{\partial
x_1^{p+1}}(\tau)d\tau\biggl |^2 dx_1\\
&+\frac{1}{2}\int_0^{\frac{1}{2\sqrt{\lambda}}}\bigg|\frac{\partial^{p}u}{\partial
x_1^{p}}\bigg|^2dx_1 +2\int_0^{\frac{1}{2\sqrt{\lambda}}}\biggl
|\int_{x_1^0}^{x_1}\frac{\partial^{p+1}u}{\partial
x_1^{p+1}}(\tau)d\tau\biggl |^2dx_1\\
= &\frac{3}{2}\int_{0}^{\frac{1}{2\sqrt{\lambda}}} \bigg
|\frac{\partial^{p}u}{\partial
x_1^{p}}\bigg|^2dx_1+3\int_0^{\frac{1}{2\sqrt{\lambda}}}\biggl
|\int_{x_1^0}^{x_1}\frac{\partial^{p+1}u}{\partial
x_1^{p+1}}(\tau)d\tau\biggl |^2 dx_1\\
\leq &\frac{3}{2}\int_{0}^{\frac{1}{2\sqrt{\lambda}}}
\bigg|\frac{\partial^{p}u}{\partial
x_1^{p}}\bigg|^2dx_1+3\int_0^{\frac{1}{2\sqrt{\lambda}}}
(x_1-x_1^0)\int_{x_1^0}^{x_1} \bigg|\frac{\partial^{p+1}u}{\partial
x_1^{p+1}}(\tau)\bigg|^2d\tau dx_1\\
\leq &\frac{3}{2}\bigg\|\frac{\partial^{p}u}{\partial
x_1^{p}}\bigg\|^2_{L^2([0,\frac{1}{2\sqrt{\lambda}}])}+3\bigg\|\frac{\partial^{p+1}u}{\partial
x_1^{p+1}}\bigg\|^2_{L^2([0,\frac{1}{2\sqrt{\lambda}}])}\int_0^{\frac{1}{2\sqrt{\lambda}}}
|x_1-x_1^0|dx_1\\
\leq&\frac{3}{2}\biggl[\bigg\|\frac{\partial^{p}u}{\partial
x_1^{p}}\bigg\|^2_{L^2([0,\frac{1}{2\sqrt{\lambda}}])}+\frac{1}{4\lambda}\bigg\|\frac{\partial^{p+1}u}{\partial
x_1^{p+1}}\bigg\|^2_{L^2([0,\frac{1}{2\sqrt{\lambda}}])}\biggl].
\endaligned$$
Thus,
\begin{equation}\label{maxu}
\aligned \underset{{x_1}}{\rm
max}\biggl|\frac{\partial^{p}u}{\partial x_1^{p}}\biggl|&\leq
\sqrt{\frac{3}{2}\biggl(2\sqrt{\lambda}\biggl\|\frac{\partial^{p}u}{\partial
x_1^{p}}\biggl\|^2_{L^2([0,\frac{1}{2\sqrt{\lambda}}])}+\frac{1}{2\sqrt{\lambda}}\biggl\|\frac{\partial^{p+1}u}{\partial
x_1^{p+1}}\biggl\|^2_{L^2([0,\frac{1}{2\sqrt{\lambda}}])}\biggl)}\\
&=\sqrt{3\biggl(\sqrt{\lambda}\biggl\|\frac{\partial^{p}u}{\partial
x_1^{p}}\biggl\|^2_{L^2([0,\frac{1}{2\sqrt{\lambda}}])}+\frac{1}{4\sqrt{\lambda}}\biggl\|\frac{\partial^{p+1}u}{\partial
x_1^{p+1}}\biggl\|^2_{L^2([0,\frac{1}{2\sqrt{\lambda}}])}\biggl)}\\
&\leq\sqrt{3\biggl(\sqrt{\lambda}\cdot2\beta_p^2\lambda^{p+n-\frac{1}{2}}
+\frac{1}{4\sqrt{\lambda}}\cdot2\beta_{p+1}^2\lambda^{p+1+n-\frac{1}{2}}\biggl)}\\
&=\sqrt{3}\cdot\lambda^{\frac{p+n}{2}}\sqrt{2\beta_p^2
+\frac{1}{2}\beta_{p+1}^2}\\
&\leq \sqrt{6}\lambda^{\frac{p+n}{2}}\sqrt{\beta_p^2
+\beta_{p+1}^2}.
\endaligned
\end{equation}
Let $u(x_1)=v_1(x_1)+\sqrt{-1}~v_2(x_1)$. We consider separately the
following two cases: \\
Case 1:  If ${\rm max} |u|\geq 6$, then at least one of ${\rm
max}|v_1|$ and ${\rm max}|v_2|$ is greater than or equal to 3.
Without loss of generality, we assume that ${\rm max} |v_1|\geq 3$
and apply Lemma 2.3 to the function $v_1$. If $v_1$ satisfies
 $${\rm max}|v_1'|< 4^{p+1}\lambda^{\frac{1}{2}}{\rm max}|v_1|,$$
 then there exists a subinterval
 $[t_1,t_2]$ of the length
 $2^{-2p-3}\lambda^{-\frac12}$ on which $|v_1|\geq \frac{1}{2}{\rm max}|v_1|\geq
 \frac32$. In fact, we can choose two points $t_1$, $t_2$ of the interval
 $\bigl[0,\frac{1}{2\sqrt{\lambda}}\bigl]$ such that $|v_1(t_2)|={\rm max}|v_1|$, $|v_1(t_1)|=\frac{1}{2}{\rm
 max}|v_1|$ and $|v_1|\geq \frac{1}{2}{\rm max}|v_1|$ on
 $[t_1,t_2]$. By the mean value theorem, we get
 $$\aligned \frac{1}{2}{\rm max}|v_1|&\leq |v_1(t_2)-v_1(t_1)|\\
 &=|v_1'(\xi)|\cdot |t_2-t_1|\\
 &\leq {\rm max}|v_1'|\cdot |t_2-t_1|\\
 &\leq 4^{p+1}\lambda^{\frac12}{\rm max}|v_1|\cdot |t_2-t_1|.\endaligned$$
 Hence,
  $$ |t_2-t_1|\geq 2^{-2p-3}\lambda^{-\frac12}.$$
  Since $\bigl|u(x_1)-e^{\sqrt{-1}\xi_1 x_1}\bigl|\geq |u(x_1)|-1\geq |v_1(x_1)|-1\geq
\frac12$ on the subinterval $[t_1,t_2]$,
  \begin{equation}
  \int_0^{\frac{1}{2\sqrt{\lambda}}}\bigl|u(x_1)-e^{\sqrt{-1}\xi_1x_1}\bigl|^{2}dx_1
   \geq 2^{-2p-5}\lambda^{-\frac12}.
   \end{equation}
  On the other hand, if $v_1$ satisfies
  $$ {\rm max}|v_1'| \leq 2^{p-1}\biggl(\frac{{\rm
max}|v_1^{(p)}|}{{\rm max}|v_1|}\biggl)^{\frac{1}{p}}{\rm
max}|v_1|,$$ then it follows from the mean value theorem and
\eqref{maxu} that the length of the subinterval of
$[0,\frac{1}{2\sqrt{\lambda}}]$ on which $|v_1|\geq \frac{3}{2}$, is
at least $2^{-p}3^{\frac1p}6^{-\frac{1}{2p}}(\beta_p^2
+\beta_{p+1}^2)^{-\frac{1}{2p}}\lambda^{-\frac{1}{2}-\frac{n}{2p}}$,
which yields
  \begin{equation}
  \int_0^{\frac{1}{2\sqrt{\lambda}}}\bigl|u(x_1)-e^{\sqrt{-1}\xi_1x_1}\bigl|^{2}dx_1 \geq 2^{-p-2}3^{\frac1p}
6^{-\frac{1}{2p}}(\beta_p^2
+\beta_{p+1}^2)^{-\frac{1}{2p}}\lambda^{-\frac12-\frac{n}{2p}}.
\end{equation}
Case 2:  Assume now that ${\rm max}|u|<6$. The latter means that
 ${\rm max}|v_1|<6$ and ${\rm max}|v_2|<6$. Since $v_1(0)=v_2(0)=0,$ there
 exists a subinterval $I$ of $[0,\frac{1}{2\sqrt{\lambda}}]$ such that $|v_1|\leq\frac13, |v_2|\leq
 \frac{1}{3}$ on $I$,
which implies
$$\bigl|u(x_1)-e^{\sqrt{-1}\xi_1x_1}\bigl|^2\geq \bigl(|u(x_1)|-1\bigl)^2\geq
\biggl(1-\frac{\sqrt{2}}{3}\biggl)^2\geq \frac14.$$ Next we apply
Lemma 2.3 to the functions $v_1$ and $v_2$. We find out that
  the length of $I$ is at least
 $${\rm min}\biggl\{3^{-2}2^{-2p-3}\lambda^{-\frac12}, \ \
 3^{-2}2^{-p}6^{\frac{1}{2p}}(\beta_p^2
+\beta_{p+1}^2)^{-\frac{1}{2p}}\lambda^{-\frac12-\frac{n}{2p}}\biggl\}.$$
 This implies
 $$ \int_{0}^{\frac{1}{2\sqrt{\lambda}}}\bigl|u(x_1)-e^{\sqrt{-1}\xi_1x_1}\bigl|^{2}dx_1 \geq \frac{\lambda^{-\frac{1}{2}}}{9}{\rm
min}\biggl\{2^{-2p-5}, \ \
 2^{-p-2}6^{\frac{1}{2p}}(\beta_p^2
+\beta_{p+1}^2)^{-\frac{1}{2p}}\lambda^{-\frac{n}{2p}}\biggl\}. $$
This completes the proof of the lemma.\endproof

{\it Proof of Theorem 1:} \ \ According to Lemma 2.4, for each $l_i$
and $p$, we have
$$\biggl\|u-e^{\sqrt{-1}<\xi, x>}\biggl\|^2_{L^{2}(T_{l_j})}  \geq \frac{1}{9\cdot 2^{n-1}}{\rm
min}\biggl \{2^{-2p-5}\lambda^{-\frac n2}, \ 2^{-p-2
}6^{\frac{1}{2p}}(\beta_p^2+\beta_{p+1}^2)^{-\frac{1}{2p}}\lambda^{-\frac
n2-\frac{n}{2p}} \biggl\},$$ where
$$\beta_{p+1}^2=\biggl(\frac{n+2}{4n\pi^2}\biggl)^{\frac{n}{2}}B_nV^{2}
D_p.$$ Next we estimate the sequence $D_{p}.$ A direct inspection
shows that
$$D_p< 2^{(2n+18)p^2}.$$
This implies that
$$
(\beta_p^2+\beta_{p+1}^2)^{-\frac{1}{2p}}\geq
(2\beta_{p+1}^2)^{-\frac{1}{2p}}
>2^{-\frac{1}{2p}-(n+9)p}\biggl(\frac{4n\pi^2}{n+2}\biggl)^{\frac{n}{4p}}B_n^{-\frac{1}{2p}}V^{-\frac1p}.
$$ Hence
$$\biggl\|u-e^{\sqrt{-1}<\xi, x>}\biggl\|^2_{L^{2}(T_{l_i})}  \geq \frac{1}{9\cdot 2^{n-1}}{\rm
min}\biggl \{2^{-2p-5}\lambda^{-\frac n2}, \ 2^{-2-(n+10)p
}\bigg(\frac{V\lambda^{\frac
n2}}{\alpha_1}\bigg)^{-\frac{1}{p}}\lambda^{-\frac n2} \biggl\},$$
where
$$ \alpha_1=\sqrt{\frac{3}{B_n}\left(\frac{4n\pi^2}{n+2}\right)^\frac
n2}.$$ Now we choose
$$p=\left[\sqrt{\frac{{\rm log}_{2}\big[(V/\alpha_1)^{n-1}
\lambda^{\frac n2}\big]}{n+12}}~\right],$$ then
$$2^{-2-(n+10)p
}\geq\bigg(\frac{V}{\alpha_1}\bigg)^{-\frac{n-2}{2p}}\bigg(\frac{V\lambda}{\alpha_1}\bigg)^{-\frac{n}{2p}}.$$
Therefore, we obtain
\begin{equation}\label{thin}
\biggl\|u-e^{\sqrt{-1}<\xi, x>}\biggl\|^2_{L^{2}(T_{l_i})} \geq
\frac{1}{9\cdot 2^{n-1}}\biggl \{
\biggl(\frac{V}{\alpha_1}\biggl)^{\frac{n}{2}}\biggl(\frac{V\lambda}{\alpha_1}\biggl)^{-\frac{n}{2}-\frac
np} \biggl\}.
\end{equation}
 Notice that for each $i$ the number of
these $n$-dimensional rectangles $T_{l_i}$ is at least
$$N_i=\biggl[\frac{A_i}{6}\lambda^{\frac{n-1}{2}}\biggl].$$
Summing the inequality \eqref{thin} for all
$l_i=1,\cdots,N_i$ and all $i=1,\cdots,m$, we immediately get
$$\aligned
V-(2\pi)^n F_\lambda(\xi)&= \biggl\|u-e^{\sqrt{-1}<\xi,
x>}\biggl\|^2_{L^{2}(\Omega)}\\
 &\geq \frac{1}{9\cdot 2^{n-1}}\biggl \{
\biggl(\frac{V}{\alpha_1}\biggl)^{\frac{n}{2}}\biggl(\frac{V\lambda}{\alpha_1}\biggl)^{-\frac{n}{2}-\frac
np}
\biggl\}\sum_{i=1}^m\biggl[\frac{A_i}{6}\lambda^{\frac{n-1}{2}}\biggl]
\Theta
(\lambda-\lambda_0)\\
&\geq \frac{\alpha_1^{-\frac12}V^{\frac12}}{9^2\cdot 2^{n}}
\biggl(\frac{V\lambda}{\alpha_1}\biggl)^{-\frac{1}{2}-\frac
np}A(\partial\Omega) \Theta(\lambda-\lambda_0),\endaligned$$ where
$$\lambda_0={\rm max}\left\{\frac{4n}{\underset{i}{\rm min}\{d_i^2\}},
\ \left(\frac{\alpha_1}{V}\right)^{\frac{2}{n}}, \
2^{\frac{2(n+12)}{n}}\left(\frac{\alpha_1}{V}\right)^{\frac{2(n-1)}{n}},
\ \left(\frac{12}{\underset{i}{\rm min}
\{A_i\}}\right)^{\frac{2}{n-1}}\right\}.$$ This yields the following
upper bound on $F_\lambda$:
\begin{equation}F_\lambda(\xi)\leq
\frac{V}{(2\pi)^n}\biggl[1-\alpha_2V^{-\frac12}A(\partial\Omega)\biggl(\frac{V\lambda}{\alpha_1}\biggl)^{-\frac{1}{2}-\frac
np} \Theta(\lambda-\lambda_0)\biggl]\triangleq M(\lambda),
\end{equation}
where
$$\alpha_2=\frac{\alpha_1^{-\frac12}}{9^2\cdot 2^{n}}.$$
By Lemma 2.1, we obtain
$$
\aligned &\sum_{j=1}^{N(\lambda)}
 \lambda_j
 =\int_{\Omega}|\xi|^2F_\lambda(\xi)d\xi\\
 \geq &\frac{n}{n+2}B_n^{-\frac2n}N(\lambda)^{\frac{n+2}{n}}M(\lambda)^{-\frac2n}\\
=&\frac{n}{n+2}B_n^{-\frac2n}N(\lambda)^{\frac{n+2}{n}}(2\pi)^2V^{-\frac2n}
\biggl[1-\alpha_2V^{-\frac12}A(\partial\Omega)\biggl(\frac{V\lambda}{\alpha_1}\biggl)^{-\frac{1}{2}-\frac
np}\Theta(\lambda-\lambda_0)\biggl]^{-\frac2n}\\
\geq&\frac{n}{n+2}B_n^{-\frac2n}N(\lambda)^{\frac{n+2}{n}}(2\pi)^2V^{-\frac2n}
\biggl[1+\frac{2\alpha_2}{n}V^{-\frac12}A(\partial\Omega)\biggl(\frac{V\lambda}{\alpha_1}\biggl)^{-\frac{1}{2}-\frac
np}\Theta(\lambda-\lambda_0)\biggl]\\
=&\frac{n}{n+2}\frac{(2\pi)^2}{(B_nV)^{\frac2n}}N(\lambda)^{\frac{n+2}{n}}+\frac{2\alpha_2}{n+2}\frac{(2\pi)^2}{(B_nV)^{\frac2n}}
N(\lambda)^{\frac{n+2}{n}}V^{-\frac12}A(\partial\Omega)\biggl(\frac{V\lambda}{\alpha_1}\biggl)^{-\frac{1}{2}-\frac
np}\Theta(\lambda-\lambda_0).
\endaligned$$

Taking $k$ large enough such that $\lambda=\lambda_k$, we
immediately get
$$\frac1k\sum_{j=1}^{k} \lambda_j
\geq\frac{n}{n+2}\frac{(2\pi)^2}{(B_nV)^{\frac2n}}k^{\frac{2}{n}}
+\frac{8\pi^2}{9^22^n(n+2)B_n^{\frac2n}} \frac{ \
A(\partial\Omega)}{V^{1+\frac
2n}}k^{\frac{2}{n}}\lambda_k^{-\frac{1}{2}}\biggl(\frac{V\lambda_k}{\alpha_1}\biggl)^{-n\varepsilon(k)}\Theta
\biggl(\lambda_k-\lambda_0\biggl),
$$ where
$$ \ \ \ \ \ \ \ \ \ \ \ \ \ \ \ \ \ \ \ \ \ \ \ \ \ \ \ \ \ \ \ \
\varepsilon(k)=\frac{1}{p} =\left[\sqrt{\frac{{\rm
log}_{2}\big[(V/\alpha_1)^{n-1} \lambda_k^{\frac
n2}\big]}{n+12}}~\right]^{-1}. \ \ \ \ \ \ \ \ \ \ \ \ \ \ \ \ \ \ \
\ \ \ \ \ \ \ \ \ \ \ \ \ \ \ \ \ \ \endproof$$


\providecommand{\bysame}{\leavevmode\hbox
to3em{\hrulefill}\thinspace}

\end{document}